\documentclass{amsart}%
\usepackage{amsthm,amstext,amscd,latexsym,amsmath,amsfonts,amssymb}
\usepackage{amsmath}
\usepackage{amsfonts}
\usepackage{amssymb}

\usepackage{graphicx}%
\theoremstyle{plain}
\newtheorem{theorem}{Theorem}
\newtheorem{lemma}[theorem]{Lemma}
\newtheorem{proposition}[theorem]{Proposition}

\theoremstyle{definition}

\newtheorem{corollary}[theorem]{Corollary}
\newtheorem{example}{Example}
\theoremstyle{remark}

\begin{document}

\title{Around the Hossz\'u-Gluskin theorem for ${n}$-ary groups}
\author{Wies{\l }aw A. Dudek}
\address{Institute of Mathematics and Computer Science,
Wroc{\l}aw University of Techno\-logy, Wybrze\.ze Wyspia\'nskiego
27, 50-370 Wroc{\l}aw, Poland} \email{dudek@im.pwr.wroc.pl}
\author{\fbox{\;Kazimierz G{\l}azek\;}}
\address{Faculty of Mathematics, Computer Science and Econometrics,
University of Zielona G\'ora, ul. Podg\'orna 50, Zielona G\'ora,
Poland} \email{K.Glazek@wmie.uz.zgora.pl} \maketitle \sloppy
\maketitle 

\centerline{\it Dedicated to Boles{\l}aw Gleichgewicht on his 85th birthday}

\begin{abstract}
We survey results related to the important Hossz\'u-Gluskin
Theorem on $n$-ary groups adding also several new results and
comments. The aim of this paper is to write all such results in
uniform and compressive forms. Therefore some proofs of new
results are only sketched or omitted if their completing seems to
be not too difficult for readers. In particular, we show as the
Hossz\'u-Gluskin Theorem can be used for evaluation how many
different $n$-ary groups (up to isomorphism) exist on some small
sets. Moreover, we sketch as the mentioned theorem can be also
used for investigation of $\mathcal{Q}$-independent subsets of
semiabelian $n$-ary groups for some special families $\mathcal{Q}$
of mappings. Such investigations will be continued.
\end{abstract}

\section{Introduction}

The non-empty set $G$ together with an $n$-ary operation $f:G^n\to
G$ is called an {\it $n$-ary groupoid} (or an {\it $n$-ary
operative}  -- in the Gluskin terminology, cf. \cite{glu65}) and
is denoted by $(G;f)$. We will assume that $n>2$.

According to the general convention similar to that introduced in
the theory of $n$-ary systems by G. \v{C}upona (cf. \cite{Cup})
the sequence of elements $x_i,x_{i+1},\ldots,x_j$ is denoted by
$x_i^j$. In the case $j<i$ it is the empty symbol. If
$x_{i+1}=x_{i+2}=\ldots=x_{i+t}=x$, then instead of
$x_{i+1}^{i+t}$ we write $\stackrel{(t)}{x}$. In this convention
$f(x_1,\ldots,x_n)= f(x_1^n)$ and
 \[
 f(x_1,\ldots,x_i,\underbrace{x,\ldots,x}_{t},x_{i+t+1},\ldots,x_n)=
 f(x_1^i,\stackrel{(t)}{x},x_{i+t+1}^n) .
 \]
If $\,m=k(n-1)+1$, then the m-ary operation $\,g\,$ of the form
 \[
 g(x_{1}^{k(n-1)+1}):=\underbrace{f(f(...,f(f}_{k}
(x_{1}^{n}),x_{n+1}^{2n-1}),...),x_{(k-1)(n-1)+2}^{k(n-1)+1})
\]
is denoted by $\,f_{(k)}$ and is called the {\it simple iteration}
of the operation $f$ (cf. \cite{Mar'64}) or an $m$-ary operation
{\it derived} from $f$. In certain situations, when the arity of
$\,g\,$ does not play a crucial role or when it will differ
depending on additional assumptions, we write $\,f_{(.)}\,$ to
mean $\,f_{(k)}\,$ for some $\,k=1,2,...$.

An $n$-ary groupoid $(G;f)$ is called {\it $(i,j)$-associative} if
 \begin{equation}
f(x_1^{i-1},f(x_i^{n+i-1}),x_{n+i}^{2n-1})=
f(x_1^{j-1},f(x_j^{n+j-1}),x_{n+j}^{2n-1})\label{assolaw}
 \end{equation}
holds for all $x_1,\ldots,x_{2n-1}\in G$. If this identity holds
for all $1\leqslant i<j\leqslant n$, then we say that the
operation $f$ is {\it associative} and $(G;f)$ is called an {\it
$n$-ary semigroup} (or, in the Gluskin's terminology, an {\it
$n$-ary associative}). It is clear that an $n$-ary groupoid is
associative iff it is $(1,j)$-associative for all $j=2,\ldots,n$.
In the binary case (i.e. for $n=2$) it is a usual semigroup.

If for all $x_0,x_1,\ldots,x_n\in G$ and fixed
$i\in\{1,\ldots,n\}$ there exists an element $z\in G$ such that
\begin{equation}                                                          \label{solv}
f(x_1^{i-1},z,x_{i+1}^n)=x_0 ,
\end{equation}
then we say that this equation is {\it $i$-solvable} or {\it
solvable at the place $i$}. If this solution is unique, then we
say that (\ref{solv}) is {\it uniquely $i$-solvable}.

An $n$-ary groupoid $(G;f)$ uniquely solvable for all
$i=1,\ldots,n$ is called an {\it $n$-ary quasigroup}. An
associative $n$-ary quasigroup is called an {\it $n$-ary group}.
For an $n$-ary quasigroup with the non-empty center to be an
$n$-ary group it is sufficient to postulate the
$(i,j)$-associativity for some fixed $1\leqslant i<j\leqslant n$
(cf. \cite{Du'86}). It is clear that for $n=2$ we obtain a usual
group.

Note by the way that in many papers $n$-ary semigroups ($n$-ary
groups) are called $n$-semigroups ($n$-groups, respectively).
Moreover, in many papers, where the arity of the basic operation
does not play a crucial role, we can find the term a {\it polyadic
semigroup} ({\it polyadic group}) (cf. \cite{Busse}, \cite{post},
\cite{Sokh}).

Now such and similar $n$-ary systems have many applications in
different branches. For example, in the theory of automata
\cite{Busse} $n$-ary semigroups and $n$-ary groups are used, some
$n$-ary groupoids are applied in the theory of quantum groups
\cite{Nik}. Different applications of ternary structures in
physics are described by R. Kerner in \cite{Ker}. In physics there
are used also such structures as $n$-ary Filippov algebras (see
\cite{Poj}) and $n$-Lie algebras (see \cite{Vai}). Some $n$-ary
structures induced by hypercubes have application in
error-correcting and error-detecting coding theory, cryptology, as
well as in the theory of $(t,m,s)$-nets (see for example
\cite{LaMu} and \cite{LaMuWhi}).

The idea of investigations of such groups seems to be going back
to E. Kasner's lecture \cite{kasner} at the fifty-third annual
meeting of the American Association for the Advancement of Science
in 1904. But the first paper concerning the theory of $n$-ary
groups was written (under inspiration of Emmy Noether) by W.
D\"ornte in 1928 (see \cite{dor}). In this paper D\"ornte observed
that any $n$-ary groupoid $(G;f)$ of the form
 $\, f(x_1^n)=x_1\circ x_2\circ\ldots\circ x_n$, where $(G;\circ
 )$ is a group, is an $n$-ary group but for every $n>2$ there are
$n$-ary groups which are not of this form. $n$-ary groups of the
first form are called {\it reducible} to the group $(G;\circ )$ or
{\it derived} from the group $(G;\circ )$, the second one are
called {\it irreducible}. Moreover, in some $n$-ary groups there
exists an element $e$ (called an {\it $n$-ary identity} or {\it
neutral element}) such that
 \begin{equation}                                            \label{n-id}
 f(\stackrel{(i-1)}{e},x,\stackrel{(n-i)}{e})=x
 \end{equation}
holds for all $x\in G$ and for all $i=1,\ldots,n$. It is
interesting that $n$-ary groups containing a neutral element are
reducible (cf. \cite{dor}). Irreducible $n$-ary groups do not
contain such elements. On the other hand, there are $n$-ary groups
with two, three and more neutral elements. The set
$\mathbb{Z}_{n-1}=\{0,1,\ldots,n-2\}$ with the operation
$f(x_1^{n})=(x_1+x_2+\ldots +x_{n})({\rm mod}\,(n-1))$ is a simple
example of an $n$-ary group in which every element is neutral. All
$n$-ary groups with this property are derived from the commutative
group of the exponent $k|(n-1)$.

It is worthwhile to note that in the definition of an $n$-ary
group, under the assumption of the associativity of $f$, it
suffices only to postulate the existence of a solution of
(\ref{solv}) at the places $i=1$ and $i=n$ or at one place $i$
other than $1$ and $n$. Then one can prove the uniqueness of the
solution of (\ref{solv}) for all $i=1,\ldots,n$ (cf. \cite{post},
p. $213^{17}$).

The above definition of $n$-ary groups is a generalization of the
Weber's and Huntington formulation of axioms of a group as a
semigroup in which the equations $xa=b$, $ya=b$ have solutions.
Many authors used the notion of $n$-ary groups as a generalization
of Pierpont's definition of groups as a semigroup with neutral and
inverse elements. Unfortunately, in this case we obtain only
$n$-ary groups derived from groups.

E.I. Sokolov proved in \cite{sok} that in the case of $n$-ary
quasigroups (i.e. in the case of the existence of a unique
solution of (\ref{solv}) at any place $i=1,\ldots,n$) it is
sufficient to postulate the $(j,j+1)$-associativity for some fixed
$j=1,\ldots,n-1$.

Using the same method as Sokolov we can prove the following
proposition (for details see \cite{DGG}):

\begin{proposition}\label{DGG1}
An $n$-ary groupoid $(G;f)$ is an $n$-ary group if and only if
$($at least$)$ one of the following conditions is satisfied:
\begin{enumerate}
\item [$(a)$] the $(1,2)$-associative law holds and the
equation $(\ref{solv})$ is solvable for $\,i=n\,$ and uniquely
solvable for $\,i=1$,
\item [$(b)$] the $(n-1,n)$-associative law holds and the
equation $(\ref{solv})$ is solvable for $\,i=1\,$ and uniquely
solvable for $\,i=n$,
\item [$(c)$] the $(i,i+1)$-associative law holds for some
$\,i\in \{2,...,n-2\}\,$ and the equation $(\ref{solv})$ is
uniquely solvable for $\,i\,$ and some $j>i$.
\end{enumerate}
\end{proposition}

\medskip

In \cite{DD80} (see also \cite{cel77}) the following
characterization of $n$-ary groups is given:

\begin{proposition}
An $n$-ary semigroup $(G;f)$ is an $n$-ary group if and only if
for some $k\in\{1,2,\ldots,n-2\}$ and all $a_1^k\in G$ there are
elements $x_{k+1}^{n-1},\,y_{k+1}^{n-1}\in G\,$ such that
\begin{equation}
f(a_1^k,x_{k+1}^{n-1},b)=f(b,y_{k+1}^{n-1},a_1^k)=b\label{DG-80}
\end{equation}
for all $\,b\in G$.
\end{proposition}

\medskip

\begin{proposition}
An $n$-ary semigroup $(G;f)$ is an $n$-ary group if and only if
for some $i,j\in\{1,2,\ldots,n-1\}$ and all $a,b\in G$ there are
$x,y\in G\,$ such that
\begin{equation}
f(x,\stackrel{(i-1)}{b},\stackrel{(n-i)}{a})=f(\stackrel{(n-j)}{a},\stackrel{(j-1)}{b},y)=b.
\label{gal-r1}
\end{equation}
\end{proposition}

\medskip
Putting in the above proposition $i=j=1$ we obtain the following
main result of \cite{Tyu85}.

\begin{corollary}
{\em An $n$-ary semigroup $(G;f)$ is an $n$-ary group if and only
if for all $a,b\in G$ there are $x,y\in G\,$ such that}
\[
f(x,\stackrel{(n-1)}{a})=f(\stackrel{(n-1)}{a},y)=b.
\]
\end{corollary}

\bigskip

From the definition of an $n$-ary group $(G;f)$ we can directly
see that for every $x\in G$ there exists only one $z\in G$
satisfying the equation
 \begin{equation}                                               \label{skew}
f(\stackrel{(n-1)}{x},z)=x .
 \end{equation}
This element is called {\it skew} to $x$ and is denoted by
$\overline{x}$. In a ternary group ($n=3$) derived from the binary
group $(G;\circ)$ the skew element coincides with the inverse
element in $(G;\circ)$. Thus, in some sense, the skew element is a
generalization of the inverse element in binary groups. This
suggests that for $n\geqslant 3$ any $n$-ary group $(G;f)$ can be
considered as an algebra $(G;f,\bar{\,}\;)$ with two operations:
one $n$-ary $\,f:G^n\to G$ and one unary $\;\bar{\,} :
x\to\overline{x}$. D\"ornte proved (see \cite{dor}) that in
ternary groups for all $x\in G$ we have
$\overline{\overline{x}}=x$, but for $n>3$ this is not true. For
$n>3$ there are $n$-ary groups in which one fixed element is skew
to all elements (cf. \cite{D90}) and $n$-ary groups in which any
element is skew to itself. Then, in the second case, of course the
$n$-ary group operation $f$ is idempotent. An $n$-ary group in
which $f(\stackrel{(n)}{x})=x$ for every $x\in G$ is called {\it
idempotent}.

Nevertheless, the concept of skew elements plays a crucial role in
the theory of $n$-ary groups. Namely, as D\"ornte proved, the
following theorem is true.

\begin{theorem}\label{dor-th}
 In any $n$-ary group $(G;f)$ the following identities:
 \begin{eqnarray}
  f(\stackrel{(i-2)}{x},\overline{x},\stackrel{(n-i)}{x},y)=y, \label{dor-r}\\
  f(y,\stackrel{(n-j)}{x},\overline{x},\stackrel{(j-2)}{x})=y, \label{dor-l}\\
  f(\stackrel{(k-1)}{x},\overline{x},\stackrel{(n-k)}{x})=x
  \label{skew2}
  \end{eqnarray}
hold for all $\,x,y\in G$, $\,2\leqslant i,j\leqslant n\,$ and
$\,1\leqslant k\leqslant n$.
\end{theorem}

\medskip

The first two identities, called now {\it D\"ornte's identities},
are used by many authors to describe the class of $n$-ary groups.
For example, in 1967 B. Gleichgewicht and K. G{\l}azek proved in
\cite{GG67} (see also \cite{Sio67}) that for fixed $n\geqslant 3$
the class of all $n$-ary groups, considered as algebras of type
$(n,1)$, forms a Mal'cev variety and found the system of
identities defining this variety. This means that all congruences
of a given $n$-ary group commute and that the lattice of all
congruences of a fixed $n$-ary group is modular. But, as was
observed many years later, from the theorem on page 448 in
Gluskin's paper \cite{glu65} it follows that the system of
identities given by B. Gleichgewicht and K. G\l azek is not
independent. For similar axiom considerations, see also
\cite{cel77}, \cite{rus79} and \cite{rus81} (for other systems of
axioms, see, e.g., \cite{monk2}). The first independent system of
identities defining this variety was given in our paper
\cite{DGG}. Now we give the minimal system of such identities.
This is the main result of \cite{Rem}.

\begin{theorem}\label{DGG}
The class of $n$-ary groups coincides with the variety of $n$-ary
groupoids $(G;f,\bar{}\;)$ with a unary operation
$\,\bar{}:x\to\overline{x}$ satisfying for some fixed
$i,j\in\{2,\ldots,n\}$ the D\"ornte identities $(\ref{dor-r})$,
$(\ref{dor-l})$ and the identity
\[
 f(f(x_1^{n}),x_{n+1}^{2n-1})=f(x_1,f(x_2^{n+1}),x_{n+2}^{2n-1}).
 \]
\end{theorem}

\medskip

Theorem \ref{DGG} gives the minimal system of identities defining
$n$-ary groups. In fact, for $n>3$ the set $Z$ of all integers
with the operation $f(x_1^n)=x_{n-1}+x_n$ is an example of a
$(1,2)$-associative $n$-ary groupoid in which (\ref{dor-r}) holds
for $\overline{x}=0$ but (\ref{dor-l}) is not satisfied.
Similarly, $(Z;f)$ with $f(x_1^n)=x_1$ is an example of a
$(1,2)$-associative $n$-ary groupoid satisfying (\ref{dor-l}) but
not (\ref{dor-r}). It is clear that the $(1,2)$-associativity
cannot be deleted.

Note by the way that in some papers there are investigated
so-called {\it infinitary} semigroups and quasigroups, i.e.
groupoids $(G;f)$, where the number of variables in the operation
$f:G^{\infty}\to G$ is infinite, but countable. Infinitary
semigroups are the infinitary groupoids $(G;f)$, where for all
natural $\,i, j\,$ the operation $f$ satisfies the identity
 \[
f(x_1^{i-1},f(x_i^{\infty}),y_1^{\infty})=
f(x_1^{j-1},f(x_j^{\infty}),y_1^{\infty}).
 \]
Infinitary quasigroups are infinitary groupoids $(G;f)$ in which
the equation $\,f(x_1^{k-1},z_k,x_{k+1}^{\infty})=x_0\,$ has a
unique solution $z_k$ at any place $k$.

From the general results obtained in \cite{belz} and \cite{MTC}
one can deduce that infinitary groups have only one element. Below
we present a simple proof of this fact.

If $(G;f)$ is an infinitary group, then, according to the
definition, for any $y,z\in G$ and $u=f(\stackrel{(\infty)}{y})$
there exists $x\in G$ such that
$z=f(u,y,x,\stackrel{(\infty)}{y})$. Thus
\[\arraycolsep.5mm
 \begin{array}{rl}
f(z,\stackrel{(\infty)}{y})=&f(f(u,y,x,\stackrel{(\infty)}{y}),\stackrel{(\infty)}{y})=
 f(u,y,f(x,\stackrel{(\infty)}{y}),\stackrel{(\infty)}{y})\\[6pt]
=&f(f(\stackrel{(\infty)}{y}),y,f(x,\stackrel{(\infty)}{y}),\stackrel{(\infty)}{y})
=f(y,f(\stackrel{(\infty)}{y}),y,f(x,\stackrel{(\infty)}{y}),\stackrel{(\infty)}{y})\\[6pt]
=&f(y,u,y,f(x,\stackrel{(\infty)}{y}),\stackrel{(\infty)}{y})
=f(y,f(u,y,x,\stackrel{(\infty)}{y}),\stackrel{(\infty)}{y})=f(y,z,\stackrel{(\infty)}{y}),
 \end{array}
 \]
i.e. for all $y,z\in G\,$ we have
\[
f(z,\stackrel{(\infty)}{y})=f(y,z,\stackrel{(\infty)}{y}).
\]
Using this identity and the fact that for all $\,x,y\in G$ there
exists $z\in G$ such that $\,x=f(z,\stackrel{(\infty)}{y}),\,$ we
obtain
\[\arraycolsep.5mm
 \begin{array}{rl}
f(\stackrel{(\infty)}{x})=&f(x,f(z,\stackrel{(\infty)}{y}),\stackrel{(\infty)}{x})=
 f(x,f(y,z,\stackrel{(\infty)}{y}),\stackrel{(\infty)}{x})\\[6pt]
=&f(x,y,f(z,\stackrel{(\infty)}{y}),\stackrel{(\infty)}{x})
=f(x,y,\stackrel{(\infty)}{x}),
 \end{array}
 \]
which together with the existence of only one solution at the
second place implies $x=y$. Hence $G$ has only one element.

\medskip

According to Theorem~\ref{DGG}, the class of all $n$-groups (for
$n>2$) can be considered as a variety of algebras
$(G;f,\bar{\,}\;)$ with one $n$-ary operation $f$ and one unary
$\,\bar{}:x\to\overline{x}$. The class of $n$-ary groups can be
also considered as a variety of algebras of different types (cf.
\cite{filomat} and \cite{jus'03}).

Theorem~\ref{DGG} is valid for $n>2$, but, as it was observed in
\cite{Rem}, this theorem can be extended to the case $n=2$.
Namely, let \ $\hat{\,}:x\to\hat{x}$ be a unary operation, where
$\hat{x}$ is defined as a solution of the equation
$f_{(2)}(\stackrel{(2n-2)}{x},\hat{x})=x$. Then using the same
method as in the proof of Theorem 2 in \cite{DGG} we can prove:

\begin{theorem}
Let $(G;f)$ be an $n$-ary $(n\geqslant 2)$ semigroup with a unary
operation $\hat{}:x\to\hat{x}$. Then $(G;f,\,\hat{\,}\,)$ is an
$n$-ary group if and only if for some $i,j\in\{2,\ldots,2n-1\}$
the following identities
\[
f_{(2)}(y,\stackrel{(i-2)}{x},\hat{x},\stackrel{(2n-1-i)}{x})=y=
f_{(2)}(\stackrel{(2n-1-j)}{x},\hat{x},\stackrel{(j-2)}{x},y)
\]
hold.
\end{theorem}

From this theorem we can deduce other definitions of $n$-ary
$(n\geqslant 2)$ groups presented in \cite{Gal'95}, \cite{Gal'03},
\cite{rus79} and \cite{Sio67}.

\medskip

An $n$-ary group is said to be {\em semiabelian} if the following
identity
\begin{equation}
f(x_1^n)=f(x_n,x_2^{n-1},x_1)
\end{equation}
is satisfied. In this case the operation $\,\bar{\,}:x\to
\overline{x}\,$ is a homomorphism (cf. \cite{GG'77}). Note by the
way that the class of all semiabelian $n$-ary groups coincides
with the class of {\it medial} $n$-ary groups (cf. \cite{medial},
[29]). (Some authors used also the name {\it abelian} instead of
{\it semiabelian} (see, e.g., \cite{Sio67}, [29]).) Such $n$-ary
groups are a special case of {\it $\sigma$-permutable} $n$-ary
groups (cf. \cite{St-D'86}), i.e. $n$-ary groups in which
$f(x_1^n)=f(x_{\sigma(1)},x_{\sigma(2)},\ldots,x_{\sigma(n)})$ for
fixed $\sigma\in S_n$. An $n$-ary group which is
$\sigma$-permutable for every $\sigma\in S_n$ is usually called
{\it commutative}.

An {\it $n$-ary power} of $x$ in an $n$-ary group $(G;f)$ is
defined in the following way: $x^{<0>}=x$, \
$x^{<1>}=f(\stackrel{(n)}{x})\;$ and
$x^{<k+1>}=f(\stackrel{(n-1)}{x},x^{<k>})\;$ for all $k>0$. In
this convention $x^{<-k>}$ means an element $z$ such that
$f(x^{<k-1>},\stackrel{(n-2)}{x},z)=x^{<0>}=x$. Then
$\overline{x}=x^{<-1>}$ and
\[
\begin{array}{l}
f(x^{<k_1>},\ldots,x^{<k_n>})=x^{<k_1+\ldots+k_n+1>}\\[4pt]
(x^{<k>})^{<t>}=x^{<kt(n-1)+k+t>}
\end{array}
\]
(cf. \cite{post}, \cite{Gl'82} or \cite{auto}).

Now, putting $x^{\!\!\!\!-(0)}=x$ and denoting by
$x^{\!\!\!\!-(m+1)}$ the skew element to $x^{\!\!\!\!-(m)}$, from
the above two identities and results obtained by W. A. Dudek in
\cite{auto} and \cite{medial} we deduce the following proposition.
\begin{proposition}
In any $n$-ary group \ $x^{\!\!\!\!-(m)}=x^{<S_m>}$, where
$S_m=\frac{(2-n)^m-1}{n-1}$.
\end{proposition}

This means that for every $n>2$ we have
$\overline{\overline{x}}=x^{<n-3>}$. In particular,
$\overline{\overline{x}}=x^{<1>}$ in all $4$-ary groups, and
$\overline{\overline{x}}=x^{<2>}$ in all $5$-ary groups (cf. \cite{Gl'82}).

\section{Hossz\'u-Gluskin algebras}

Let $(G;f)$ be an $n$-ary group. Fixing in $f(x_1^n)$ some $m<n$
elements we obtain a new $(n-m)$-ary operation which in general is
not associative. It is associative only in the case when these
fixed elements are located in some special places, for example, in
the case when this new operation has the form
\begin{equation}
g(x_1^k)=f(x_1,a_1^r,x_2,a_1^r,x_3,a_1^r,\ldots,a_1^r,x_{k-1},a_1^r,x_k),
\label{k-ary}
 \end{equation}
where $k+r(k-1)=n$ and $a_1,\ldots,a_r\in G$ are fixed. Of course,
in this case $(G;g)$ is a $k$-ary group. It is denoted by
$ret_{a_1^r}(G;f)$ and is called a {\it $k$-ary retract of}
$(G;f)$ (see, e.g., \cite{DM2}). 
 For different elements
$a_1,\ldots,a_r$ we obtain different $k$-ary retracts, but all
these $k$-ary retracts (for a fixed $n$-ary group) are isomorphic
(cf. \cite{DM2}). 
Therefore, we can consider only retracts for
$a_1=\ldots=a_r=a$. In such retracts the element skew to $x$ has
the form
\[
f_{(\cdot)}(\stackrel{(n-r-2)}{a},\overline{a},\stackrel{(n-3)}{x},\overline{x},
\stackrel{(n-r-2)}{a},\overline{a},\stackrel{(n-3)}{x},\overline{x},\ldots,
\stackrel{(n-r-2)}{a},\overline{a},\stackrel{(n-3)}{x},\overline{x},
\stackrel{(n-r-2)}{a},\overline{a}),
\]
where $\overline{a}$ and $\overline{x}$ are skew in $(G;f)$. This
means that the skew elements in this retract can be expressed by
the operations of an $n$-ary group $(G;f)$.

A very important role play {\it binary retracts}, especially
retracts denoted by $ret_a(G;f)$, where $x\circ
y=f(x,\stackrel{(n-2)}{a},y)$. The identity of the group
$(G;\circ)$ is $\overline{a}$. One can verify that the inverse
element to $x$ has the form
\begin{equation}
x^{-1}=f(\overline{a},\stackrel{(n-3)}{x},\overline{x},\overline{a}).\label{inv}
\end{equation}
Thus, in this group
 \begin{equation}
x\circ
y^{-1}=f(x,\stackrel{(n-3)}{y},\overline{y},\overline{a}).\label{x-y}
 \end{equation}

Binary retracts of an $n$-ary group $(G;f)$ are commutative only
in the case when $(G;f)$ is semiabelian (medial). So, $(G;f)$ is
semiabelian if and only if
 \[
f(x,\stackrel{(n-2)}{a},y)=f(y,\stackrel{(n-2)}{a},x)
 \]
holds for all $x,y\in G$ and some fixed $a\in G$.

M. Hossz\'u was first who observed a strong connection between
$n$-ary groups and their binary retracts. He proved in
\cite{Hosszu} the following theorem:

\begin{theorem}\label{thGH}
An $n$-ary groupoid $(G;f)$, $n>2$, is an $n$-ary group if and
only if
\begin{enumerate}
\item[$(i)$] on $G$ one can define a binary operation $\cdot$ such
that $(G;\cdot)$ is a group,
\item[$(ii)$] there exist an automorphism $\varphi$ of $(G;\cdot)$ and
$b\in G$ such that $\varphi(b)=b$,
\item[$(iii)$] $\varphi^{n-1}(x)=b\cdot x\cdot b^{-1}$ holds for every
$x\in G$,
\item[$(iv)$] $f(x_1^n)=x_1\cdot\varphi(x_2)\cdot\varphi^2(x_3)\cdot\varphi^3(x_4)\cdot\ldots
\cdot\varphi^{n-1}(x_n)\cdot b$ for all $x_1,\ldots,x_n\in G$.
\end{enumerate}
\end{theorem}

Two years later, this theorem was proved by L. M. Gluskin (see
\cite{glu65}) in a more general form (for so-called positional
operatives). For a generalization to $n$-ary semigroups, see also
\cite{Monk} and \cite{Zup}.  In another version this theorem was
also formulated by E. L. Post (cf. \cite{post}, p. 246). An
elegant short proof was given by E. I. Sokolov in \cite{sok}. His
proof is based on the observation that $(G;\cdot)=ret_a(G;f)$.
Then we have:
  \begin{equation}
\varphi(x)=f(\overline{a},x,\stackrel{(n-2)}{a})\label{autom}
 \end{equation}
and
 \begin{equation}
 b=f(\stackrel{(n)}{\overline{a}}).\label{b}
 \end{equation}
From (14) and (7) or (8), we can deduce that commutative $n$-ary group
operations have the form $f(x_1^n)=x_1\cdot x_2\cdot\ldots\cdot x_n\cdot b$,
where $(G;\cdot)$ is a commutative group.

Note that the last condition of Theorem~\ref{thGH} can be
rewritten in the form
\begin{equation}
f(x_1^n)=x_1\cdot\varphi(x_2)\cdot\varphi^2(x_3)\cdot\varphi^3(x_4)\cdot\ldots
\cdot\varphi^{n-2}(x_{n-1})\cdot b\cdot x_n. \label{e10}
\end{equation}

The above theorem has the following generalization (cf.
\cite{DM1}):
\begin{theorem}\label{genGH}
An $n$-ary groupoid $(G;f)$, $n>2$, is an $n$-ary group if and
only if
\begin{enumerate}
\item[$(i)$] on $G$ one can define a $k$-ary operation $g$ such
that $(G;g)$ is a $k$-ary group and $k-1$ divides $n-1$,
\item[$(ii)$] there exist an automorphism $\varphi$ of $(G;g)$
and elements $b_2,\ldots,b_k\in G$ such that $\varphi(b_i)=b_i$
for $i=2,\ldots,k$,
\item[$(iii)$] $g(\varphi^{n-1}(x),b_2^k)=g(b_2^k,x)$ holds for every
$x\in G$,
\item[$(iv)$]
$f(x_1^n)=g_{(\cdot)}(x_1,\varphi(x_2),\varphi^2(x_3),\ldots,
\varphi^{n-1}(x_n),b_2^k)$ for all $x_1,\ldots,x_n\in G$.
\end{enumerate}
\end{theorem}

In this theorem $(G;g)=ret_{a_1^r}(G;f)$, where $a_1=\ldots=a_r=a$.
In this case, we get:
\begin{equation}
\varphi(x)=f(\overline{a},\stackrel{(n-r-2)}{a},x,\stackrel{(r)}{a}),
\label{k-fi}
\end{equation}
\begin{equation}
b_2=
f_{(\cdot)}(\stackrel{(n-r-2)n}{a},\stackrel{(n)}{\overline{a}},\stackrel{(k-2)(n-r-2)}{a}),
 \ \ b_3=\ldots=b_k=\overline{a}.\label{b_2}
\end{equation}

Other important generalizations can be found in \cite{Hosszu2}
(for heaps), \cite{Mar-Jan} (for vector valued groups),
\cite{Sokh} (for partially associative $n$-ary quasigroups).

\medskip

Following E. L. Post (see \cite{post}, cf. [4], p. 36--40, and
[28]) a binary group \ $\mathfrak{G}^{\ast }=\left(G^{\ast };\circ
\right)$ is said to be a {\it covering group} for the $n$-ary
group $(G;f)$ if there exists an embedding $\tau :G\rightarrow
G^{\ast }$ such that $\tau(G)$ is a generating set of $G^{\ast }$
and \ $\tau (f(x_1^n)) = \tau
(x_1)\circ\tau(x_2)\circ\ldots\circ\tau(x_n)$ for every
$x_1,\ldots,x_n\in G.$ \ $\mathfrak{G}^{\ast }$ is a
\textit{universal covering group} (or a \textit{free covering
group}) if for any covering group $\mathfrak{G}_{1}^{\ast }$\
there exists a homomorphism from $G^*$ onto $G^*_1$ such that the
following diagram is commutative (or compatible -- in another
terminology):

\begin{center}
\begin{minipage}{6cm}

\ \ \ \ \ \ \ \ \ \ \ \ \ \ \ \ \ \ \ \ $G$

\ \ \ \ \ \ \ \ \ \ \ \ \ \ {  \ }$\swarrow ${  \ }$
\circlearrowleft ${  \ }$\searrow $

\ \ \ \ \ \ \ \ \ \ \ \ \ $G^{\ast}--\rightarrow \ G_{1}^{\ast}$

\ \ \ \ \ \ \ \ \ \ \ \ \ \ \ \ \ \ \ \ {\footnotesize onto}
\end{minipage}
\end{center}

Post proved in \cite{post} that for every $n$-ary group $(G;f)$
there exist a covering group $(G^{\ast};\circ)$ and its normal
subgroup $G_0$ such that $G^{\ast}\diagup G_0$ is a cyclic group
of order $n-1$ and $f(x_1^n)=x_1\circ x_2\circ\ldots\circ x_n$ for
all $x_1,\ldots,x_n\in G$. So, the theory of $n$-ary groups is
closely related to the theory of {\it cyclic extensions of
groups}, but these theories are not equivalent.

Indeed, the above theorems show that for any $n$-ary group $(G;f)$
we have the sequence
$$
O\rightarrow (G_0;\circ)\rightarrow (G^{\ast};\circ )
\stackrel{\zeta }{\longrightarrow }\mathcal{C}\left( n\right)
\rightarrow O,
 $$
where $(G^{\ast };\circ)$ is the free covering group of $(G;f)$
with $G=\zeta^{-1}\left( 1\right)$, and $1$ is a generator of the
cyclic (additively writing) group $\mathcal{C}\left( n\right)
=(C_{n};+_{n}).$

We have
\[
\begin{array}{ccccc}
&  & ( G_{1}^{\ast };\circ) &  &  \\
& \nearrow & \uparrow & \searrow &  \\
(G_0;\circ ) &  & \circlearrowright \;\;\; {\vert}
\circlearrowleft\circlearrowright &  & C\left( n\right) \\
& \searrow &\downarrow & \nearrow &  \\
&  & (G_2^{\ast};\circ)&  &
\end{array}
\]
where we use

$\circlearrowright $ \ for the equivalence of extensions,

$\circlearrowleft $ \ for the isomorphism of suitable $n$-ary
groups.

\medskip

Of course, two $n$-ary groups determined in the above-mentioned
sense by two equivalent cyclic extensions are isomorphic. However,
two non-equivalent cyclic extensions can determine two isomorphic
$n$-ary groups.

\begin{example} Consider two cyclic extensions of the cyclic group $%
\mathcal{C}(3)$ by $\mathcal{C}(3)$:
$$
0\rightarrow \mathcal{C}(3)\stackrel{\alpha }{\rightarrow
}\mathcal{C}(9)\stackrel{\beta _{1}}{\rightarrow
}\mathcal{C}(3)\rightarrow 0
$$
and
$$
0\rightarrow \mathcal{C}(3)\stackrel{\alpha }{\rightarrow
}\mathcal{C}(9)\stackrel{\beta _{2}}{\rightarrow
}\mathcal{C}(3)\rightarrow 0,
$$
where the homomorphisms $\alpha $, $\beta _{1}$ and $\beta _{2}$
are given by:
\[
\begin{array}{lcl}
\;\alpha (x)=3x&{\rm for}& x\in C_{3},\\[4pt]
\beta_{1}(x)\equiv x({\rm mod}\,3)&{\rm for}&x\in C_{9},\\[4pt]
\beta_{2}(x)\equiv 2x({\rm mod}\,3)&{\rm for}&x\in C_{9}.
\end{array}
\]

It is easy to verify that if $g(x,y,z,v)=(x+y+z+v)({\rm mod}\,9),$
then the $4$-ary groups $(\beta _{1}^{-1}(1);g)$ and $(\beta
_{2}^{-1}(2);g)$, corresponding to those extensions are isomorphic
to the $4$-ary groups $(C_{3};f_{1})$, $(C_{3};f_{2})$,
respectively, where
 $$
f_{1}(x,y,z,v)\equiv (x+y+z+v+1)({\rm mod}\,3)
 $$
and
$$
f_{2}(x,y,z,v)\equiv (x+y+z+v+2)({\rm mod}\,3).
$$
These 4-groups are isomorphic. The isomorphism
$\varphi:(C_{3};f_{1})\to(C_{3};f_{2})$ has the form
$\varphi(x)\equiv 2x({\rm mod}\,3)$. Nevertheless, the
above-mentioned extensions are not equivalent (because there is no
automorphism $\lambda $ of $\mathcal{C}(9)$ such that $\lambda
\circ \alpha =\alpha $ and \ $\beta _{2}\circ \lambda =\beta
_{1}$).
\end{example}

The algebra $(G;\cdot,\varphi,b)$ of the type $(2,1,0)$, where
$(G;\cdot)$ is a (binary) group, $b\in G$ is fixed, $\varphi\in
Aut(G;\cdot)$, $\varphi(b)=b$ and $\varphi^{n-1}(x)=b\cdot x\cdot
b^{-1}$ for every $x\in G$ is called a {\it Hossz\'u-Gluskin
algebra} (briefly: an {\it $HG$-algebra}). We say that an
$HG$-algebra $(G;\cdot,\varphi,b)$ is {\it associated} with an
$n$-group $(G;f)$ if the identity (\ref{e10}) is satisfied. In
this case we say also that an $n$-ary group $(G;f)$ is {\it
$\langle\varphi,b\rangle$-derived} from the group $(G;\cdot)$. A
$k$-ary $HG$-algebra $(G;g,\varphi,b_2^k)$ can be defined
similarly. Binary $HG$-algebras are studied in \cite{jus'95a},
\cite{jus'95b} and \cite{jus'03}.

Theorems~\ref{thGH} and \ref{genGH} state that every $k$-ary
$HG$-algebra is associated with some $n$-ary group. Any $n$-ary
group is $\langle\varphi,b\rangle$-derived from some binary group
and $\langle\varphi,b_2^k\rangle$-derived from some $k$-ary group.

\section{Calculation of $n$-ary groups on small sets}

Results presented in the previous section give the possibility to
evaluate the number of non-isomorphic $n$-ary groups. To calculate
these groups we must use the following result proved in
\cite{DM2}.

\begin{theorem}\label{izoth}
Two $n$-ary groups $(G_1;f_1)$, $(G_2;f_2)$ are isomorphic if and
only if for every $c\in G_1$ there exists an isomorphism
$h:ret_c(G_1;f_1)\to ret_d(G_2;f_2)$ such that $d=h(c)$, \
$h(f_1({\overline{c}},\ldots,{\overline{c}}))=
f_2(\overline{d},\ldots,\overline{d}) $ and
$h(f_1(\overline{c},x,\!\stackrel{(n-2)}{c}))=
f_2(\overline{d},h(x),\!\stackrel{(n-2)}{d}).$
\end{theorem}

\begin{corollary}\label{izoth2}
{\it Two commutative $n$-ary groups $(G_1;f_1)$, $(G_2;f_2)$ are
isomorphic if and only if for every $c\in G_1$ there exists an
isomorphism $h:ret_c(G_1;f_1)\to ret_d(G_2;f_2)$ such that
$d=h(c)$ and $h(f_1({\overline{c}},\ldots,{\overline{c}}))=
f_2(\overline{d},\ldots,\overline{d})$.}
\end{corollary}

If $(G,\cdot)$ is an abelian group, then, of course, we can
consider the automorphism of the form $\alpha(x)=x^{-1}$. Then $G$
with the operation
\begin{equation}
f(x_1^n)=x_1\cdot x_2^{-1}\cdot x_3\cdot\ldots\cdot
x^{-1}_{n-1}\cdot x_n\label{-1}
\end{equation}
is an $n$-ary group if $n$ is odd. Such $n$-ary groups are
characterized by the following theorem proved in \cite{Gl-Mi'84}.

\begin{theorem} Let $m$ be odd and let $(G;f)$ be an
$n$-ary group. Then the operation $f$ has the form $(\ref{-1})$,
where $(G;\cdot)$ is an abelian group, if and only if
\begin{enumerate}
\item[$(i)$] \ $f\left( x,\ldots,x\right) =x$,
\item[$(ii)$] \ $f\left(x_1^{i},y,y,x_{i+3}^{n}\right)=
f\left(x_1^{i},z,z,x_{i+3}^{n}\right)$ for all $\;0\leqslant
i\leqslant n-2.$
\end{enumerate}
In this case $(G;\cdot)=ret_a(G;f)$ for some $a\in G$.
\end{theorem}

Consider an abelian group $(G;+)$. Then, as a special case of an $n$-ary
group operation of form (19), one can obtain the ternary term operation
$$
f(x,y,z)=x-y+z
$$
which is a so-called Mal'tsev term in the group $G$. Of course, it is
idempotent and medial ({\it entropic} -- in another terminology). Such
ternary operations appear in several branches of mathematics. For example,
they play very important role in affine geometry and the theory of modes
(because of idempotency and mediality), in the theory of congruences
in general algebras (because existence of a Mal'tsev term in general algebras
implies permutability of congruences and then modularity of lattices of
congruences) and also in the theory of clones which is important in
Universal Algebra and as well in Multiple-valued Logics.

\medskip
From results obtained in \cite{Gl-Mi'84} (cf. also \cite{Tim'72})
we can deduce:

\begin{proposition} Let $(G;\cdot)$ be a group and let $t_1,\ldots,t_n$
be fixed integers. Then $G$ with the operation
\[
f(x_1^{m})=(x_1)^{t_{1}}\cdot(x_{2})^{t_{2}}\cdot\ldots \cdot
(x_{n-1})^{t_{n-1}}\cdot(x_{n})^{t_{n}},
\]
is an $n$-ary group if and only if
\begin{enumerate}
\item \ $x^{t_1}=x=x^{t_n}$,
\item \ $t_j=k^{j}$ \ for some integer $k$ and all $j=2,\ldots,n-1$,
\item \ $(x\cdot y)^k=x^k\cdot y^k.$
\end{enumerate}
\end{proposition}

In this case we say that $(G;f)$ is {\em derived from the
$k$-exponential group}.
\begin{proposition}
An $n$-ary group $(G;f)$ is derived from the $k$-exponential
$(k>0)$ group $(G;\cdot)$ if and only if there exists $a\in G$
such that
\begin{enumerate}
\item \ $f(a,\ldots,a)=a$,
\item \ $f_{(k)}(\stackrel{(n-2)}{a},x,\stackrel{(n-2)}{a},x,
\ldots,\stackrel{(n-2)}{a},x,a)=x$.
\end{enumerate}
Moreover, $(G;\cdot)=ret_a(G;f)$.
\end{proposition}

Using the above results we can describe all non-isomorphic $n$-ary
groups with small numbers of elements.

For this let $(\mathbb{Z}_k;+)$ be the cyclic group modulo $k$.
Consider the following $n$-ary operation:
 \[\arraycolsep=.5mm
 \begin{array}{rl}
f_{a}( x_1^n) &\equiv (x_{1}+\ldots+x_{n}+a)\,({\rm mod}\,k) ,\\[4pt]
g_{d}(x_1^n)&\equiv (x_1+dx_2+\ldots+d^{n-2}x_{n-1}+x_n)\,({\rm mod}\,k),\\[4pt]
g_{d,c}(x_1^n)&\equiv
(x_1+dx_2+\ldots+d^{n-2}x_{n-1}+x_n+c)\,({\rm mod}\,k),
\end{array}
 \]
where $a\in\mathbb{Z}_k$, \ $c,d\in\mathbb{Z}_k\setminus\{0,1\},$
\ $d^{\,n-1}\equiv 1\,({\rm mod}\,k).$ Additionally, for the
operation $g_{d,c}$ we assume that $dc\equiv c\,({\rm mod}\,k)$
holds. By Theorem~\ref{thGH}, $(\mathbb{Z}_{k};f_{a})$,
$(\mathbb{Z}_{k};g_{d})$ and $(\mathbb{Z}_{k};g_{d,c})$ are
$n$-ary groups.

\bigskip

In \cite{Gl-Mi-S} the following theorem is proved:
\begin{theorem} A $k$-element $n$-ary group $(G;f)$ is
$\langle\varphi,b\rangle$-derived from the cyclic group of order
$k$ if and only if it is isomorphic to exactly one $n$-ary group
of the form $(\mathbb{Z}_{k};f_{a})$, $(\mathbb{Z}_{k};g_{d})$ or
$(\mathbb{Z}_{k};g_{d,c})$, where $d|gcd(k,n-1)$ and $c|k.$
\end{theorem}

An infinite cyclic group can be identified with the group
$(\mathbb{Z};+)$. This group has only two automorphisms:
$\varphi(x)=x$ and $\varphi(x)=-x$. So, according to
Theorem~\ref{thGH}, $n$-ary groups defined on $\mathbb{Z}$ have
the form $(\mathbb{Z};f_a)$ or $(\mathbb{Z};g_{-1})$, where
$$
g_{-1}(x_1^n)=x_1-x_2+x_3-x_4+\ldots+x_n, $$ and $n$ is odd. Since
$\varphi_k(x)=x+k$ is an isomorphism of $n$-ary groups
$(\mathbb{Z};f_a)$ and $(\mathbb{Z};f_b)$, where $a=b+(n-1)k$, the
calculation of non-isomorphic $n$-ary groups of the form
$(\mathbb{Z};f_a)$ can be reduced to the case when
$a=0,1,\ldots,n-2$. From Corollary~\ref{izoth2} it follows that
these $n$-ary groups are non-isomorphic.

So, we have proved
\begin{theorem}
An $n$-ary group $\langle\varphi,b\rangle$-derived from the
infinite cyclic group $(\mathbb{Z};+)$ is isomorphic to an $n$-ary
group $(\mathbb{Z};f_a)$, where $0\leqslant a\leqslant (n-2)$, or
to $(\mathbb{Z};g_{-1})$, where $n$ is odd.
\end{theorem}

Denote by $Inn\left(G;\cdot\right)$ the group of all inner
automorphisms of $(G;\cdot)$, by $Out\left( G;\cdot \right)$ the
factor group of $Aut\left(G;\cdot\right)$ by
$Inn\left(G;\cdot\right)$, and by $Out_{n}\left( G;\cdot \right) $
the set of all cosets $\overline{\gamma}\in
Out\left(G;\cdot\right)$ containing $\gamma$ such that
$\gamma^{n-1}\in Inn\left(G;\cdot\right)$. Then, as it is proved
in \cite{Gl-Mi-S}, for centerless groups, i.e. groups for which
$\,card(Cent\left(G;\cdot\right))=1$, the following theorem is
true.

\begin{theorem} Let $\left( G;\cdot \right)$ be a centerless group such that
$Out_{n}\left(G;\cdot\right)$ is abelian, and let
$\left(G;f\right)$ be $\langle\alpha,a\rangle$-derived, and
$\left(G;g\right)$ be $\langle\beta,b\rangle$-derived from
$\left(G;\cdot\right)$. Then $\left(G;f\right)$ is isomorphic to
$\left( G;g\right)$ if and only if $\,\alpha\circ\beta^{-1}\in
Inn\left( G;\cdot\right)$.
\end{theorem}

The number of pairwise non-isomorphic $n$-ary groups
$\langle\varphi,b\rangle$-derived from a centerless group
$(G;\cdot)$ is smaller or equal to $s=card(
Out_{n}\left(G;\cdot\right))$. It is equal to $s$ if and only if
$Out\left(G;\cdot\right)$ is abelian.

For every $n$ and $k\neq 2,6$, there exists exactly one $n$-ary
group which is $\langle\varphi,b\rangle$-derived from $S_{k}$ (for
$k=2$ and $k=6$ we have one or two such $n$-ary groups relatively
to evenness of $n$).

Let now $(G;\cdot)$ be an arbitrary group, $c\in G,$ $\varphi \in
Aut\left(G;\cdot \right) $. Let us put
\[
\arraycolsep=.5mm \begin{array}{rl}
 f_{c}^{(\cdot)}(x_1^n)&=x_1\cdot x_2\cdot\ldots\cdot x_n\cdot
 c,\\[4pt]
g_{\varphi }^{(\cdot)}(x_1^n)& =x_1\cdot\varphi\left(x_{1}\right)
\cdot\ldots\cdot\varphi^{n-1}( x_n) ,\\[4pt]
g_{\varphi ,c}^{( \cdot)}(x_1^n)& =x_1\cdot\varphi\left(
x_{2}\right)\cdot\ldots\cdot\varphi^{n-1}(x_n)\cdot c.
\end{array}
\]

For example (for details see \cite{Gl-Mi'87}), we have the
following:

\begin{theorem} Let $l=gcd\left(n-1,12\right),$ \ $\left(G_{4};\ast\right)$
be the Klein four-group $($with $0$ as the neutral element$)$, let
$\gamma ,\varepsilon \in Aut\left( G_{4};\ast \right),$ where
$\gamma $ is of order $2$ and $\varepsilon $ of order $3$, and let
$c\in G_{4}\backslash \{0\}$ be the fix point of $\gamma $. Then
every $n$-ary group $\langle\varphi,b\rangle$-derived from
$(G_{4};\ast)$ is isomorphic to exactly one $(G_{4};f)$, where $f$
is one of the following $n$-ary group operations:
\begin{enumerate}
\item[$(a)$] \ $f_{0}^{\left( \ast \right) },$ $f^{\left( \ast
\right) },g_{\ \gamma }^{\left( \ast \right) },g_{\ \gamma
,c}^{\left( \ast \right) } $ or $g_{\ \varepsilon }^{\left( \ast
\right) }$ \ \ for \ $l=12$,
\item[$(b)$] \ $f_{0}^{\left( \ast \right) },$ $f_{1}^{\left(
\ast \right) },g_{\ \gamma }^{\left( \ast \right) }$ or $g_{\
\varepsilon }^{\left( \ast \right) }$ \ \ for \ $l=6$,
\item[$(c)$] \ $f_{0}^{\left( \ast \right) },$ $f_{1}^{\left( \ast
\right) },g_{\ \gamma }^{\left( \ast \right) }$ or $g_{\ \gamma
,c}^{\left( \ast \right) }$ \ \ for \ $l=4$,
\item[$(d)$] \ $f_{0}^{\left( \ast \right) }$ or $g_{\ \varepsilon
}^{\left( \ast \right) }$ \ \ for \ $l=3$,
\item[$(e)$] \ $f_{0}^{\left( \ast \right) },$ $f_{1}^{\left( \ast
\right) }$ or $g_{\ \gamma }^{\left( \ast \right) }$ \ \ for \
$l=2$, \item[$(f)$] \ $f_{0}^{\left( \ast \right) }$ \ \ for \
$l=1.$
\end{enumerate}
\end{theorem}

Comparing our results with results obtained in \cite{Gl-Mi'87},
\cite{Gl-Mi'88} and \cite{Gl-Mi-S} (cf. also \cite{post} for
$k=2,3$) we can tabularize the numbers of $n$-ary groups on
$k$-element sets with $k<8$ in the following way (we use the
abbreviations: commut. = commutative, idem. = idempotent):

\bigskip {\small\noindent
\begin{tabular}{|l|c|c|}
\hline $k=2$, \ $l=gcd\,(n-1,2)$ & $l=2$& $l=1\rule{0mm}{3mm}$
\\ \hline
$n\equiv t({\rm mod}\,2)$&$t=1$&$t=0\rule{0mm}{3mm}$\\
\hline all &2& 1\\ \hline
 commutative&2 & 1 \\ \hline
 commutative, idempotent & 1 & 0 \\ \hline
\end{tabular}

\bigskip\noindent
\begin{tabular}{|l|c|c|c|c|}
\hline $k=3$, \ $l=gcd\,(n-1,6)$ &$l=6$ &$l=3$&$l=2$&$l=1\rule{0mm}{3mm}$ \\
\hline $n\equiv t\,({\rm mod}\,6)$&$t=1$
&$t=4$&$t=3,\,5$&$t=0,\,2\rule{0mm}{3mm}$\\ \hline
 {all} &3&2&2&1\\ \hline
{commutative} &2&2&1&1 \\ \hline
 commutative, idempotent &1&1&0&0\\ \hline
non-commut., medial, idempotent &1&0&1&0 \\ \hline
\end{tabular}

\bigskip\noindent
\begin{tabular}{|l|c|c|c|c|c|c|}
\hline $k=4$, \ $l=gcd\,(n-1,12)$ &$l=12$&$l=6$&$l=4$&$l=3$&$l=2$&$l=1\rule{0mm}{3mm}$ \\
\hline $n\equiv t\,({\rm mod}\,12)$ & $t=1$
&$t=7$&$t=5,\,9$&$t=4,\,10$&$t=3,11$&$t=t_0\rule{0mm}{3mm}$\\
\hline {all} &10&8&9&3&7&2 \\ \hline
 commutative &5&4&5&2&4&2 \\ \hline
 {commutative, idempotent} &2&1&2&0&1&0 \\ \hline
 non-commut., medial, idem. &3&2&1&1&1&0 \\ \hline
 {non-commut., medial}, non-idem., &2&2&3&0&2&0 \\ \hline
\end{tabular}

\smallskip $t_0=0,\,2,\,6,\,8$.

\bigskip\noindent
\begin{tabular}{|l|c|c|c|c|c|c|}
\hline $k=5$, \ $l=gcd\,(n-1,20)$ &$l=20$&$l=10$&$l=5$&$l=4$&$l=2$
&$l=1\rule{0mm}{3mm}$
\\ \hline $n\equiv t\,({\rm mod}\, 20)$&$t=1$ &$t=11$&$t=6,16$&$t=t_1$
& $t=t_2$ & $t=t_3\rule{0mm}{3mm}$\\ \hline
 {all}&5&3&2&4&2&1 \\ \hline
 {commutative}&2&2&2&1&1&1 \\ \hline
 {commutative, idempotent} &1&1&1&0&0&0 \\ \hline
 non-commut., idem., medial &3&1&0&3&1&0 \\ \hline
 non-commut., non-idem., medial &0&0&0&0&0&0 \\ \hline
\end{tabular}

\smallskip $t_1=5,9,13,17$, \ \

$t_2=3,7,15,19$, \ \

$t_3=0,2,4,8,10,12,14,18$.

\bigskip\noindent
\begin{tabular}{|l|c|c|c|c|}
\hline $k=6$, \ $l=gcd\,(n-1,6)$
&$l=6$&$l=3$&$l=2$&$l=1\rule{0mm}{3mm}$
\\ \hline
$n\equiv t\,({\rm mod}\,6)$&$t=1$
&$t=4$&$t=3$&$t=0,\,2\rule{0mm}{3mm}$
\\ \hline all&7&3&5&2\\ \hline
commutative&4&2&2&1 \\
\hline {commutative, idempotent}&1&0&0&0\\ \hline {medial,
idempotent, non-commut.}&1&0&1&0 \\ \hline
 non-commut., medial, non-idem.,&1&0&1&0 \\ \hline
non-medial &1&1&1&1 \\
\hline
\end{tabular}

\bigskip\noindent
\begin{tabular}{|l|c|c|c|c|c|c|c|c|}\hline
$k=7$, \ $l=gcd\,(n-1,42)$ &$l=42$&$l=21$&$l=14$&$l=7$ &$l=6$ &
$l=3$&$l=2$ &$l=1\rule{0mm}{3mm}$ \\
\hline $n\equiv t\,({\rm mod}\,42)$&$t=1$&$t=22$&$t=t_4$ &$t=t_5$ & $t=t_6$
&$t=t_7$&$t=t_8$&$t=t_9\rule{0mm}{3mm}$ \\
\hline all&7&4&3&2&6&3&2&1 \\
\hline {commutative}&2&2&2&2&1&1&1&1\\ \hline
non-com., medial, idem.,& 5&2&1&0&5&2&1&0\\
\hline commutative, idempotent&1&1&1&1&0&0&0&0\\
\hline
\end{tabular}

\smallskip

$t_4= 15, 29$,

$t_5=8, 36$,

$t_6=7,13,19,25,31,37$,

$t_7=4,10,16,28,34,40$,

$t_8= 3,5,9,11,17,21,23,33,35,39,41$,

$t_9=0,2,6,12,14,18,20,24,26,30,32,38$.
}

\section{Term equivalence of $n$-ary groups}

For any general algebra $\frak{A}=(A;\mathbb{F})$ one can define
the set $\mathbb{T}^{(n)}(\frak{A})$ of all {\it $n$-ary term
operations} as the smallest set of $n$-ary operations on $A$
containing $n$-ary projections (or $n$-ary trivial operations, in
another terminology) and closed under compositions with
fundamental operations. Then the set
$\mathbb{T}(\frak{A})=\bigcup\limits_{n=1}^{\infty}\mathbb{T}(\frak{A})$
of all {\it term operations} is the smallest set of operations on
the set $A$ containing the set $\mathbb{F}$ of fundamental
operations and all projections $e_i^{(n)}(x_1^n)=x_i$,
($i=1,2,\ldots,n$, $n=1,2,\ldots$), and closed under (direct)
compositions. Of course, $\mathbb{T}(\frak{A})$ is a {\it clone}
in the sense of Ph. Hall (see, e.g., \cite{Cohn}). It is worth
mentioning that the term operations were also called {\it
algebraic operations} by several authors (see, e.g.,
\cite{Mar'61}). Two algebras $\frak{A}_1=(A;\mathbb{F})$ and
$\frak{A}_2=(A;\mathbb{G})$ are called {\it term equivalent} if
$\mathbb{T}(\frak{A}_1)=\mathbb{T}(\frak{A}_2)$ (see, e.g.,
\cite{Gb}, p. 32, 56). If elements from some subsets $A_1$ and
$A_2$ of $A$ are treated as constant elements of algebras
$\frak{A}_1=(A;\mathbb{F}\cup A_1)$ and
$\frak{A}_2=(A;\mathbb{G}\cup A_2)$, respectively, and
$\mathbb{T}(\frak{A}_1)=\mathbb{T}(\frak{A}_2)$, then $\frak{A}_1$
and $\frak{A}_2$ are {\it polynomially equivalent}. Two varieties
$\mathcal{V}_1$ and $\mathcal{V}_2$ of algebras (perhaps of
different types) are term equivalent (polynomially equivalent,
respectively) if for every algebra $\frak{A}_1\in\mathcal{V}_1$
there exists an algebra $\frak{A}_2\in\mathcal{V}_2$ term
equivalent (polynomially equivalent, resp.) to $\frak{A}_1$, and
vice versa.

\medskip

Using Theorem~\ref{genGH} and taking into account formulas
(\ref{k-fi}) and (\ref{b_2}), we have

\begin{theorem}\label{teq}
Let $\mathfrak{G}=(G;f,\bar{}\;)$ be an $n$-ary group for a fixed
$n>2$, an element $a$ belong to $G$, and let $k$ be such a natural
number that $(k-1)$ divide $(n-1)$. Then the algebra
$\mathfrak{G}_a=(G;f,\bar{}\; ,a)$, with the additional constant
$a\in G$ is term equivalent to the algebra $(G;g,\varphi,b_2^k)$,
where $\varphi$ is an automorphism of a $k$-ary group $(G;g)$,
$(k-1)$ divides $(n-1)$, and $b_2,\ldots,b_k$ are constant
elements in $G$ such that $\varphi(b_i)=b_i$ for $i=2,\ldots,k$
and $g(\varphi^{n-1}(x),b_2^k)=g(b_2^k,x)$ for all $x\in G $.
\end{theorem}

Indeed, $f$ is determined by $g$, $\varphi$ and $b_2,\ldots,b_k$
by the formula $(iv)$ from Theorem~\ref{genGH}. The function \
$\bar{}:x\to\overline{x}$ can be easily expressed by the operation
$g$. Namely, if $f=g_{(t)}$, then $\overline{x}=x^{<-t>}$, where
$x^{<s>}$ is a $k$-ary power of $x$. According to
Theorem~\ref{genGH}, the element $\overline{x}$ also can be
expressed by $g$, $\varphi$ and $b_2,\ldots,b_n$ as a solution $z$
of the equation
$$
x=f(\stackrel{(n-1)}{x},z)=g_{(\cdot)}(x,\varphi(x),\varphi^2(x),\ldots,
\varphi(x)^{n-2},b_2^k,z).
$$

Conversely, the operations of $(G;g,\varphi,b_2^k)$ are term
derived from the operations of $(G;f,\bar{}\;)$ by (\ref{k-fi})
and (\ref{b_2}). $(G;g)=ret_{a_1^r}(G;f)$, where
$a_1=\ldots=a_r=a$, which completes the proof of
Theorem~\ref{teq}.

\medskip

By Theorem~\ref{thGH} and formulas (\ref{autom}) and (\ref{b}), we
have the following corollaries.

\begin{corollary}\label{term1}{\em
Let $\mathfrak{G}=(G;f,\bar{}\;)$ be an $n$-ary group for a fixed
$n>2$, and let an element $a$ belong to $G$. Then the algebra
$\mathfrak{G}_a=(G;f,\bar{}\; ,a)$ is term equivalent to the
$HG$-algebra $(G;\cdot,\varphi,b)$, where $(G;\cdot)$ is a group,
$\varphi\in Aut(G;\cdot)$, $b\in G$, $\varphi(b)=b$,
$\varphi^{n-1}(x)=b\cdot x\cdot b^{-1}$ for all $x\in G$.}
\end{corollary}

\begin{corollary}\label{term2}
{\em For fixed $n>2$, the variety of $n$-ary groups $($as algebras
of type $(n,1)\,)$ is polynomially equivalent to the variety of the
corresponding $HG$-algebras $($as algebras of type $(2,1,1,0)\,)$.
}
\end{corollary}

Let now $\mathfrak{G}=(G;f,\bar{}\;)$ be a semiabelian $n$-ary
group $(n>2)$. \linebreak Then the $HG$-algebra associated with
$\mathfrak{G}$ has a commutative group operation denoted by $+$.
Let $\mathfrak{H}=(G;+,\varphi,b)$ be associated with
$\mathfrak{G}$ and $\mathfrak{G}_a=(G;f,\bar{}\;,a)$. Then
$\mathfrak{H}$ and $\mathfrak{G}_a$ are term equivalent (see
Theorems \ref{thGH} and \ref{genGH},  Corollary~\ref{term1}, and
formulas (\ref{inv}) -- (\ref{b_2})). In this case we have
\[
\arraycolsep=.5mm
\begin{array}{rl}
-y&=f(\overline{a},\stackrel{(n-3)}{x},\overline{x},\overline{a}\,),\\[4pt]
x+y&=f(x,\stackrel{(n-3)}{(-y)},\overline{(-y)},\overline{a}\,),\\[4pt]
\varphi(x)&=f(\overline{a},x,\stackrel{(n-2)}{a}),\\[4pt]
{\rm and } \ \ \ \ b&=f(\stackrel{(n)}{\overline{a}}).
\end{array}
\]

We can describe of all term operations of $\mathfrak{G}_a$ by using
the language of $HG$-algebras.

At first, we consider unary term operations. Denote by $g_i(x)$
the following operation
\begin{equation}\label{g_i}
g_i(x)=k_{i1}\varphi^{l_{i1}}(x)+k_{i2}\varphi^{l_{i2}}(x)+\ldots
+k_{it}\varphi^{l_{it}}(x)
\end{equation}
for some $t,l_{i1},\ldots,l_{it}$ non-negative integers and some
$k_{i1},\ldots,k_{it}\in\mathbb{Z}$. Then it is easily to verify

\begin{lemma}
Let $\,\mathfrak{H}=(G;+,\varphi,b)$ be the $HG$-algebra
associated with a semiabelian $n$-ary group $\mathfrak{G}$. Then
all unary term operations of $\mathfrak{H}$ $($and of
$\;\mathfrak{G}_a\,)$ are of the form
\begin{equation}\label{g}
g(x)=g_i(x)+k_g b
\end{equation}
for some $g_i$ of the form $(\ref{g_i})$ and $k_g\in\mathbb{Z}$.
\end{lemma}

Indeed, it is enough to observe that
$g\in\mathbb{T}^{(1)}(\mathfrak{H})$, $\varphi(g(x))$ is again of
the form (\ref{g}), and the set of all such operations is closed
under addition.

\begin{theorem}
Let $\,\mathfrak{H}=(G;+,\varphi,b)$ be the $HG$-algebra
associated with a semiabelian $n$-ary group $\mathfrak{G}$. Then
all $m$-ary term operations of $\mathfrak{H}$ $($and of
$\;\mathfrak{G}_a\,)$ are of the form
\begin{equation}\label{sum}
F(x_1,\ldots,x_m)=\sum\limits_{i=1}^{m}g_i(x_i)+k_F b
\end{equation}
for some $g_i(x)$ of the form $(\ref{g_i})$ and
$k_F\in\mathbb{Z}$.
\end{theorem}

A verification of this theorem can be done by induction with
respect to the complexity of term operations and we left it to
readers.

\section{$\mathcal{Q}$-independent sets in $HG$-algebras}

E. Marczewski observed at the end of the 1950s that there are
common features of linear independence of vectors and
set-theoretical independence, and proposed a general scheme of
independence called here {\em $\mathcal{M}$-independence}. Recall
that the notion of set-theoretical independence (or, more
generally, independence in Boolean algebras, see, e.g.,
\cite{BalF82}, \cite{BrRu81}, \cite{Gla71}, \cite{Mar60}) was
introduced at the mid-1930s by G. Fichtenholz and L. Kantorovich
\cite{FiKa34} and also, independently, by E. Marczewski himself,
and this notion is very important in Measure Theory (see, e.g.,
\cite{FiKa34}, \cite{Mar38}, \cite{Mar48a}, \cite {Myc68},
\linebreak and \cite{Sik64}).

Let $\mathfrak{A}=(A;\mathbb{F})$ be an algebra $\emptyset \neq
X\subseteq A$. The set $X$ is said to be $\mathcal{M}$\textit{-independent}
(see \cite{Mar'58}, \cite{Mar'61}) $(X\in Ind(\mathfrak{A};\mathcal{M})$,
for short$)$ if

\begin{enumerate}
\item[(a) ]  $(\forall n\in \mathbb{N}$, $n\leqslant card(X))$
$(\forall f,g\in\mathbb{T}^{(n)}(\mathfrak{A}))$
($\forall\underset{\neq }{\underbrace{a_{1},\ldots ,a_{n}}}\in X$)

\hspace*{3mm}$\big[f(a_{1}^{n})=g(a_{1}^{n})\Longrightarrow f=g \
({\rm in }\ A)\big]$.

\bigskip

This condition is equivalent to each of the following ones
\medskip

\item[(b)]  $(\forall n\in \mathbb{N}$, $n\leqslant card(X))$
$(\forall f,g\in\mathbb{T}^{(n)}(\mathfrak{A}))$ $(\forall p:
X\rightarrow A)$ $(\forall a_{1},\ldots ,a_{n}\in X)$
$$
\big[f(a_{1}^{n})=g(a_{1}^{n})\Longrightarrow f(p(a_{1}),\ldots
,p(a_{n}))=g(p(a_{1}),\ldots ,p(a_{n}))\big],
$$

\item[(c)]  $(\forall p\in A^{X}) \ (\exists\bar{p}\in Hom(\langle
X\rangle_{\mathfrak{A}},\mathfrak{A})) \ \bar{p}|_{X}=p$, where
$\langle X\rangle_{\mathfrak{A}}$ is a subalgebra of
$\mathfrak{A}$ generated by $X$,

\item[(d)]  $\langle X\rangle _{\mathfrak{A}}$ is a $\mathbb{K}$-{\em free
algebra $\mathbb{K}$-freely generated by} $X$, where
$\mathbb{K}=\{\mathfrak{A}\}$ (or, by Birkhoff Theorem,
$\mathbb{K}=\mathcal{H}\mathcal{S}\mathcal{P}\{\mathfrak{A}\}$, a variety gene\-rated
by $\mathfrak{A}$).
\end{enumerate}

\medskip

Basic properties of $\mathcal{M}$-independence are the following
ones:
\begin{itemize}
\item  (``\textit{hereditarity}'') $X\in Ind\;(\mathfrak{A},\mathcal{M})$, \
$Y\subseteq X\Longrightarrow Y\in
Ind\;(\mathfrak{A},\mathcal{M})$,\vspace{5pt}

\item  $(\forall X\subseteq A)$ $(\forall\,{\rm finite} \ Y\subseteq
X)\;\big( Y\in Ind(\mathfrak{A},\mathcal{M})\Longrightarrow X\in
Ind(\mathfrak{A},\mathcal{M})\big)$\\
(i.e. the family $\mathbb{J}=Ind(\mathfrak{A},\mathcal{M})$ is of finite
character).
\end{itemize}

The notion of $\mathcal{M}$-independence is stronger than that of
independence with respect to the closure operator of such a kind
$X\mapsto\langle X\rangle_{\mathfrak{A}}$ (for $X\subseteq A$).

There are some notions of independence which are not special cases
of $\mathcal{M}$-independence, such as:

\noindent
\begin{tabular}{rl}
$\bullet $ & linear independence in abelian groups, \\
$\bullet $ & independence with respect to a closure operator
$\mathcal{C}$
(i.e. $\mathcal{C}$-independence), \\
$\bullet $ & stochastic independence, \\
$\bullet $ & ``independence-in-itself'' defined by J.~Schmidt (in 1962), \\
$\bullet $ & ``weak independence'' used by S.~\thinspace
\'{S}wierczkowski (in 1964).
\end{tabular}

For this reason, a general notion of independence with respect to
a family of mappings was proposed by E.~Marczewski in 1966 (and
studied in \cite{Mar'69} and \cite{Gla71}). This notion is general
enough to cover the above-mentioned kinds of independences.

Let $\emptyset\neq X\subseteq A$ and
$$
\mathcal{Q}_X\subseteq A^X=\mathcal{M}_X=\{p\;|\;p:X\rightarrow A\},
$$
$$
\mathcal{Q}(A)=\mathcal{Q}=\bigcup\{\mathcal{Q}_X \;|\;X\subseteq A \},
$$
$$
\mathcal{M}(A)=\mathcal{M}=\bigcup\{A^X \;|\;X\subseteq A\}.
$$

For an algebra $\mathfrak{A} = (A, \mathbb{F})$, a mapping \ $p:X
\rightarrow A$ belongs to $\mathcal{H}_X (\mathfrak{A})$ if and only if
there exists a homomorphism \ $\bar{p}:\langle X\rangle
_{\mathfrak{A}}\rightarrow A$ such that $\bar{p}|_{X}=p$.

The set $X$ is said to be $\mathcal{Q}$-{\it independent} ($X\in
Ind(\mathfrak{A},\mathcal{Q})$, for short) if

\begin{center}
$\mathcal{Q}_X \subseteq \mathcal{H}_X (\mathfrak{A})$
\end{center}

\noindent or, equivalently,
$$
(\forall p\in \mathcal{Q}_{X})\;(\forall\;{\rm finite }\; n\leqslant
card(X))\;(\forall f,g\in \mathbb{T}^{(n)}(\mathfrak{A}))\;
(\forall a_{1},\ldots ,a_{n}\in X)\;$$
$$\big[f(a_{1}^{n})=g(a_{1}^{n})\Longrightarrow f(p(a_{1}),\ldots
,p(a_{n}))=g(p(a_{1}),\ldots ,p(a_{n}))\big].$$

\medskip
\noindent{\bf Examples.} (In the following examples we will use a
terminology which differs from the original one.)

\begin{enumerate}
\item[1)] $\mathcal{Q}=\mathcal{M}=\bigcup \{A^{X}\;|\;X\subseteq A\};$ $\mathcal{M}$-\textit{independence}
(E.~Marczewski: \textit{general algebraic independence},
\cite{Mar'58}),

\item[2)] $\mathcal{Q}=\mathcal{G}=\bigcup \{p|_{X}\;|\;p\in A^{A}$ is \ diminishing,
$X\subseteq A\};\;\mathcal{G}$-{\it independence }(G.~Gr\"{a}tzer:
\textit{weak independence}, \cite{Grat67}), where a mapping $p$ is
called \textit{diminishing} if
$$
(\forall f,g\in \mathbb{T}^{(1)}(\mathfrak{A}))\;(\forall a\in A)
\;\big[f(a)=g(a)\Longrightarrow f(p(a))=g(p(a))\big].$$
\end{enumerate}
For abelian groups, the notion of $\mathcal{G}$-independence gives
us the well-known \textit{linear independence}.

Now we can able to obtain some results on $\mathcal{Q}$-independence (for
special families $\mathcal{Q}$ of mappings, e.g., for $\mathcal{Q}=\mathcal{M}$ and $\mathcal{G}$ ) in
$HG$-algebras of type $\mathfrak{H}=(G;+,\varphi,b)$, where $(G;+)$
is an abelian group.

In this case, the equality
\begin{equation}\label{F=G}
F_1(x_1,\ldots,x_m)=F_2(x_1,\ldots,x_m)
\end{equation}
(for two term operations of the form (\ref{sum}) in $\mathfrak{H}$)
is equivalent to the equality
\begin{equation}\label{H=0}
H(x_1,\ldots,x_m)=0,
\end{equation}
where $H\in\mathbb{T}^{(m)}(\mathfrak{H})$, i.e.
$H(x_1,\ldots,x_m)=\sum\limits_{i=1}^{m}h_i(x_i)+k_{_H} b$, and
$0$ denotes the zero of the group $(G;+)$.

Consider a subset $X$ of $G$. Let for $a_1,\ldots,a_m\in X$ the
equality
\begin{equation}\label{H_a}
H(a_1,\ldots,a_m)=0,
\end{equation}
hold. Taking into account the mapping $p:X\to\langle
X\rangle_{\mathfrak{A}}$ defined by $p(a_i)=0$ and $p(x)=x$ for
$x\in X\setminus\{a_1,\ldots,a_m\}$, we get $k_{_H}b=0$. (We
observe that such mapping $p$ belongs to families $\mathcal{M}$ and $\mathcal{G}$ .)
Therefore we have
$$
\sum\limits_{i=1}^{m} h_i(a_i)=0.
$$
Consider the mapping $q_j:X\to \langle X\rangle_{\mathfrak{A}}$
defined for fixed $j\in\{1,\ldots,m\}$ as follows:
\[
q_j(x)=\left\{\begin{array}{ccl} a_j&{\rm if }&x=a_j ,\\[4pt]
0&{\rm if }&x\ne a_j.\end{array}\right.
\]
We obtain $h_j(a_j)=0$ for all $j=1,2,\ldots,m$. (In the
considered case all $q_j$ belong to $\mathcal{M}$ and $\mathcal{G}$ .)

In particular, we can easily observe, by similar considerations,
that the following result holds:

\begin{theorem}
Let $X\subseteq G$ be a subset of the $HG$-algebra
$\mathfrak{H}=(G;+,\varphi,b)$. Then $X\in Ind(\mathfrak{H},\mathcal{G})$ if and only if
for any $m\leqslant card(X)$ for all $a_1,\ldots,a_m\in X$ and
every term operation
$H(x_1,\ldots,x_m)=\sum\limits_{i=1}^{m}h_i(x_i)+k_{_H}b$ the
equality
\begin{equation}
\sum\limits_{i=1}^{m}h_i(a_i)+k_{_H}b=0
\end{equation}
is equivalent with
\[
\left(\forall i\in\{1,\ldots,m\}\right)\,\left(h_i(a)=0\;\&\;
k_{_H}b=0\right).
\]

Moreover, $X$ is $\mathcal{M}$-independent in this $HG$-algebra
iff for all pairwise different elements $a_1,\ldots,a_m$ from $X$
equality $(26)$ implies $h_i(x)=0$ for all $i=1,2,\ldots,m$ and
$\,k_Hb=0$.
\end{theorem}


\begin{thebibliography}{90}

\bibitem{BalF82} B. Balcar and F. Fran\v{e}k, {\em Independent families
in complete Boolean algebras}, Trans. Amer. Math. Soc. {\bf 274}
(1982), $607-618.$

\bibitem{BrRu81} F.M. Brown and S. Rudeanu, {\em Consequences,
consistency, and independence in Boolean algebras}, Notre Dame J.
Formal Logic {\bf 22} (1981), $45-62.$

\bibitem{belz} V.D. Belousov and Z. Stojakovi\'c,  {\em On
infinitary quasigroups}, Publ. Inst. Math. (Beograd) {\bf 15(29)}
(1973), $31-42.$

\bibitem{b} R.H. Bruck, {\it A Survey of binary systems}, Springer-Verlag,
Berlin 1958.

\bibitem{cel77} N. Celakoski, {\em On some axioms system for
$n$-groups}, Mat. Bilten (Skopje) {\bf 1} (1977), $5-14.$

\bibitem{Cohn} P.M. Cohn, {\em Universal Algebra}, $2^{nd}$
edition, D. Reidel Publ. Co., Dordrecht 1981.

\bibitem{Cup} G. \v{C}upona, {\em Finitary associative
operations with neutral elements}, (Makedonian), Bull. Soc. Math.
Phys. Mac\'edoine {\bf 12} (1961), $15-24.$

\bibitem{dor} W. D{\"o}rnte, {\em Untersuchungen {\"u}ber einen
verallgemeinerten Gruppenbegriff}, Math. Z. {\bf 29} (1928),
$1-19.$

\bibitem{Rem} W.A. Dudek, {\em Remarks on $n$-groups}, Demonstratio
Math. {\bf 13} (1980), $165-181.$

\bibitem{auto} W.A. Dudek, {\em Autodistributive $n$-groups},
Commentationes Math. Annales Soc. Math. Polonae, Prace
Matematyczne {\bf 23} (1983), $1-11.$

\bibitem{Du'86} W.A. Dudek, {\em On $(i,j)$-associative $n$-groupoids
with the non-empty center}, Ricerche di Matematica (Napoli) {\bf
35} (1986), $105-111.$

\bibitem{medial} W.A. Dudek, {\em Medial $n$-groups and skew elements},
p. $55-80$ in ``Universal and Applied Algebra'', World Scientific,
Singapore 1989.

\bibitem{D90} W.A. Dudek, {\em On $n$-ary group with only one skew
element}, Radovi Matemati{\v{c}}ki (Sarajevo), {\bf 6} (1990), $171-175.$

\bibitem{filomat} W.A. Dudek, {\em Varieties of polyadic groups}, Filomat
{\bf 9} (1995), $657-674.$

\bibitem{DGG} W.A. Dudek, K. G{\l}azek and B. Gleichgewicht, {\em A note on
the axioms of $n$-groups}, Colloquia Math. Soc. J. Bolyai {\bf 29}
(``Universal Algebra'', Esztergom (Hungary) 1977), $195-202.$
(North-Holland, Amsterdam 1982.)

\bibitem{DD80} W.A. Dudek and I. Gro{\'z}dzi{\'n}ska, {\em On ideals in regular
$n$-semigroups}, Matemati{\v{c}}ki Bilten (Skopje) {\bf 3/4
(29/30)} (1979-1980), $35-44.$

\bibitem{DM1} W.A. Dudek and J. Michalski, {\em On a generalization of Hossz{\'u}
theorem}, Demonstratio Math. {\bf 15} (1982), $783-805.$

\bibitem{DM2} W.A. Dudek and J. Michalski, {\em  On retracts of polyadic groups},
Demonstratio Math. {\bf 17} (1984), $281-301.$

\bibitem{DM3} W.A. Dudek and J. Michalski, {\em On a generalization of a
theorem of Timm}, Demonstratio Math. {\bf 18} (1985), $869-883.$

\bibitem{Du-St'01} W.A. Dudek and Z. Stojakovi\'c, {\em On Rusakov's
$n$-ary $rs$-groups}, Czechoslovak Math. J. {\bf 51(126)} (2001), $275-183.$

\bibitem{FiKa34} G. Fichtenholz and L. Kantorovitch, {\em Sur les
op\'{e}ration dans l'espace des fonctions born\'{e}es}, Studia
Math. \textbf{5} (1934), $69-98.$

\bibitem{Gal'95} A.M. Gal'mak, {\em New axiomatics of an $n$-ary group},
(Russian), in {\em Problems in algebra and applied mathematics}, Gomel 1995,
$31-38.$

\bibitem{Gal'97} A.M. Gal'mak, {\em Post's and Gluskin-Hossz\'u's theorems},
(Russian), Gomel Univ. Press, Gomel 1997.

\bibitem{Gal'03} A.M. Gal'mak, {\em $n$-Ary groups}, (Russian), Gomel Univ.
Press, Gomel 2003.

\bibitem{Gla71} K. G\l azek, {\em Independence with respect to a family of
mappings in abstract algebras}, Dissertationes Math. \textbf{81} (1971), $1-55.$

\bibitem{Gl'82} K. G{\l}azek, {\em Remarks on polyadic groups and their weak
homomorphisms}, Contributions to General Algebra {\bf 2} (1982), $133-152.$

\bibitem{Gb} K. G{\l}azek, {\em Algebras of algebraic operations and morphisms
of algebraic systems}, {Polish}, Acta Universitatis Wratislaviensis No. {\bf 1602},
Wroclaw 1994.

\bibitem{d} K. G\l azek and B. Gleichgewicht, {\em On some method of
construction of the covering group}, (Russian), Acta Univ. Wratislaviensis
{\bf 188} (1973), $125-135$.

\bibitem{GG'77} K. G\l azek and B. Gleichgewicht, {\em Abelian $n$-groups},
Colloquia Math. Soc. J. Bolyai {\bf 29} (``Universal Algebra'',
Esztergom (Hungary) 1977), $321-329.$ (North-Holland, Amsterdam
1982.)

\bibitem{GG'85} K. G\l azek and B. Gleichgewicht, {\em On $3$-semigroups
and $3$-groups polynominal-derived from integral domains},
Semigroup Forum {\bf 32} (1985), $61-70.$

\bibitem{Gl-Mi'84} K. G\l azek and J. Michalski, {\em On polyadic groups which
are term-derived from groups},  Studia Sci. Math. Hungar. {\bf 19}
(1984), $307-315.$

\bibitem{Gl-Mi'87} K. G\l azek and J. Michalski, {\em Note on polyadic groups on
sets with at most $7$ elements},  Beitr\"age Algebra Geom. {\bf 24} (1987),
$151-158.$

\bibitem{Gl-Mi'88} K. G\l azek and J. Michalski, {\em Polyadic group
operations on small sets}, p. $85-93$ in: ``General Algebra 1988 (Krems 1988)'',
North-Holland, Amsterdam 1990.

\bibitem{Gl-Mi-S} K. G\l azek, J. Michalski and I. Sierocki, {\em On evaluation
on numbers of some polyadic groups}, Contributions to General Algebra
{\bf 3} (1984), $157-171.$

\bibitem{GG67} B. Gleichgewicht and K. G{\l}azek, {\em Remarks on $n$-groups as
abstract algebras}, Colloq. Math. {\bf 17} (1967), $209-219.$

\bibitem{glu64} L.M. Gluskin, {\em Positional operatives}, Soviet Math.,
Doklady {\bf 5} (1964) $1001-1004$ (translation from: Dokl. Akad.
Nauk SSSR {\bf 157} (1964), $767-770).$

\bibitem{glu65} L.M. Gluskin, {\em Positional operatives}, (Russian), Mat.
Sbornik (N.S.) {\bf 68} ({\bf 110}) (1965), $444-472.$

\bibitem{Grat67} G. Gr\"atzer, {\em On a new notion of independence in
universal algebra}, Colloq. Math. {\bf 17} (1967), $225-234.$

\bibitem{Busse} J.W. Grzymala-Busse, {\em Automorphisms of polyadic automata},
J. Assoc. Comput. Mach. {\bf 16} (1969), $208-219.$

\bibitem{Hosszu} M. Hossz{\'u}, {\em On the explicit form of $n$-group
operations}, Publ. Math. (Debrecen) {\bf 10} (1963), $88-92.$

\bibitem{Hosszu2} M. Hossz{\'u}, {\em The explicit form of a class of heaps},
Publ. Technical Univ. for Heavy Industry (Miskolc) {\bf 23} (1964), $265-268.$

\bibitem{kasner} E. Kasner, {\em An extension of the group concept}
(reported by L.G. Weld), Bull. Amer. Math. Soc. {\bf 10} (1904), $290-291.$

\bibitem{Ker} R. Kerner, {\em Ternary algebraic structures and their
applications in physics}, Univ. P. and M. Curie, Paris 2000.

\bibitem{LaMu} C.F. Laywine and G.L. Mullen, {\em Discrete
Mathematics Using Latin Squares}, Wiley, New York 1998.

\bibitem{LaMuWhi} C.F. Laywine, G.L. Mullen and G. Whittle,
{\em $D$-Dimensional hypercubes and the Euler and MacNeish
conjectures}, Monatsh. Math. {\bf 111} (1995), $223-238.$

\bibitem{MTC} Z. Madevski, B. Trpenovski and G. \v{C}upona G, {\em On
infinitary associative operations} (Makedonian), Bull. Soc. Math.
Phys. Mac\'edoine {\bf 15} (1964), $19-22.$

\bibitem{Mar38} E. (Szpilrajn-)Marczewski, {\em Ensemles ind\'{e}pendants et
mesures non s\'{e}parables}, C.R. Acad. Sci. Paris \textbf{207} (1938), $247-248.$

\bibitem{Mar48a} E. Marczewski, {\em Ensembles ind\'{e}pendants et leurs
applications \`{a} la th\'{e}orie de la mesure}, Fund. Math.
\textbf{35} (1948), $13-28.$

\bibitem{Mar'58} E. Marczewski, {\em A general scheme of the notios of
independence in mathematics}, Bull. Acad. Pol. Sci., Ser. Sci. Math. Astron.
Phys. {\bf 6} (1958), $731-738.$

\bibitem{Mar60} E. Marczewski, {\em Independence in algebras of sets and
Boolean algebras}, Fund. Math. {\bf 48} (1960), $135-145.$

\bibitem{Mar'61} E. Marczewski, {\em Independence and homomorphisms in abstract
algebras}, Fund. Math. {\bf 50} (1961), $45-61.$

\bibitem{Mar'64} E. Marczewski, {\em Remarks on symmetrical and
quasi-symmetrical operations}, Bull. Acad. Pol. Sci., Ser. Sci. Math. Astron.
Phys. {\bf 12} (1964), $735-737.$

\bibitem{Mar'69} E. Marczewski, {\em Independence with respect to a family of
mappings}, Colloq. Math. {\bf 22} (1969), $19-21.$

\bibitem{Mar-Jan} S. Markovski and B. Janeva, {\em Post and Hossz\'u-Gluskin
theorem for vector valued groups}, p. $77-88$ in: Proceedings of the Conference
``Algebra and Logic'' (Sarajevo 1987), Univ. Novi Sad 1989.

\bibitem{Monk}J.D. Monk and F.M. Sioson, $m$-{\it Semigroups, and function representations}, Fund. Math. {\bf 59}
 (1966), 233--241.

\bibitem{monk2} J.D. Monk and F.M. Sioson, {\it On the general theory of m-groups},
Fund. Math. {\bf 72} (1971), $233-244$.

\bibitem{Myc68} J. Mycielski, {\em Algebraic independence and measure},
Fund. Math. \textbf{61} (1968), $165-169.$

\bibitem{Nik} D. Nikshych and L. Vainerman, {\em Finite quantum groupoids and
their applications}, Univ. California, Los Angeles 2000.

\bibitem{Poj} A.P. Pojidaev, {\em Enveloping algebras of Fillipov algebras},
Comm. Algebra {\bf 31} (2003), $883-900.$

\bibitem{post} E.L. Post, {\em Polyadic groups}, Trans. Amer. Math. Soc.
{\bf 48} (1940), $208-350.$

\bibitem{rus79} S.A. Rusakov, {\em On the definition of $n$-ary groups},
(Russian), Doklady Akad. Nauk BSSR {\bf 23} (1979), $965-967.$

\bibitem{rus81} S.A. Rusakov, {\em An existence criterion for n-groups}, (Russian),
p. $77-82$ in: ``{\it Subgroup structure of finite groups $($Gomel
1977--80$)$}'', Izd. ``Nauka i Tekhnika'', Minsk 1981.

\bibitem{Schm62} J. Schmidt, {\em Einige algebraische \"Aquivalente zum
Auswahlaxiom}, Fund. Math. {\bf 50} (1962), $485-496.$

\bibitem{Sik64} R. Sikorski, {\em Boolean Algebras}, Springer-Verlag,
Berlin 1964.

\bibitem{Sio67} F.M. Sioson, {\em On free abelian $m$-groups, $I$}, Proc.
Japan. Acad. {\bf 43} (1967), $876-879.$

\bibitem{Sokh} F.M. Sokhatsky, {\em On Dudek's problems on the skew operation
in polyadic groups}, East Asian Math. J. {\bf 19} (2003), $63-71.$

\bibitem{sok} E.I. Sokolov, {\em On the Gluskin-Hossz\'u theorem for
D\"ornte $n$-groups}, (Russian), Mat. Issled. {\bf 39} (1976), $187-189.$

\bibitem{St-D'86} Z. Stojakovi\'c and W.A. Dudek, {\em On ${\sigma }$-permutable
$n$-groups}, Publ. Inst. Math. (Beograd) {\bf 40(54)} (1986), $49-55.$

\bibitem{Swi64} S. \'Swierczkowski, {\em Topologies in free algebras}, Proc.
London Math. Soc. {\bf 14} (1964), $566-576.$

\bibitem{Tim'72} J. Timm, {\em Verbandstheoretische Behandlung $n$-stelliger
Gruppen}, Abh. Math. Semin. Univ. Hamburg {\bf 37} (1972), $218-224.$

\bibitem{Tyu85} V.I. Tyutin, {\em About the axiomatics of $n$-ary groups},
(Russian), Doklady Akad. Nauk BSSR {\bf 29} (1985), $691-693.$

\bibitem{jus'95a} J. U{\v{s}}an, {\em On Hossz\'u-Gluskin algebras corresponding
to the same $n$-group}, Zb. Rad. Prirod.-Mat. Fak. Univ. u Novom
Sadu, ser. Mat. {\bf 25.1} (1995), $101-119.$

\bibitem{jus'95b} J. U{\v{s}}an, {\em Congruences of $n$-group and of associated
Hossz\'u-Gluskin algebras}, Novi Sad. J. Math. {\bf 28} (1998), $91-108.$

\bibitem{jus'03} J. U{\v{s}}an, {\em $n$-groups in the light of the neutral
operations}, Math. Moravica 2003, special issue, p. $3-162$.

\bibitem{Vai} L. Vainerman and R. Kerner, {\em On special classes of
$n$-algebras}, J. Math. Phys. {\bf 37} (1996), $2553-2565.$

\bibitem{Ziz} M. \v{Z}i\v{z}ovi\'c, {\em A topological analogue of the
Hossz\'u-Gluskin theorem}, (Serbo-Croatian), Mat. Vesnik {\bf 13(28)}
(1976), $233-239.$

\bibitem{Zup}D. Zupnik, {\it Polyadic semigroups}, Publ. Math. (Debrecen) {\bf 14} (1967), 273--279.
\end{thebibliography}
\end{document}